\documentclass[12pt]{article}  
\usepackage{amssymb,amsmath,bm,graphicx}
\usepackage{mathrsfs}
\usepackage{constants}
\usepackage{hyperref}
\usepackage{enumerate}
\usepackage{enumitem}
\usepackage{float}
\usepackage{caption}
\usepackage{verbatim}
\usepackage{tikz}
\usetikzlibrary{arrows.meta}

\providecommand{\cmark}[2][]{\relax} 

\tikzset{%
    >={Latex[width=2mm,length=2mm]},
		base/.style = {rectangle, rounded corners, draw=black,
                           minimum width=4cm, minimum height=1cm,
                           text centered, font=\sffamily}
													}

\usepackage{color}
\begin{document}
\newcommand{\bea}{\begin{eqnarray}}
\newcommand{\ena}{\end{eqnarray}}
\newcommand{\beas}{\begin{eqnarray*}}
\newcommand{\enas}{\end{eqnarray*}}
\newcommand{\beq}{\begin{equation}}
\newcommand{\enq}{\end{equation}}
\def\qed{\hfill \mbox{\rule{0.5em}{0.5em}}}
\newcommand{\bbox}{\hfill $\Box$}
\newcommand{\ignore}[1]{}
\newcommand{\ignorex}[1]{#1}
\newcommand{\wtilde}[1]{\widetilde{#1}}
\newcommand{\qmq}[1]{\quad\mbox{#1}\quad}
\newcommand{\qm}[1]{\quad\mbox{#1}}
\newcommand{\nn}{\nonumber}
\newcommand{\Bvert}{\left\vert\vphantom{\frac{1}{1}}\right.}
\newcommand{\To}{\rightarrow}
\newcommand{\E}{\mathbb{E}}
\newcommand{\Var}{\mathrm{Var}}
\newcommand{\Cov}{\mathrm{Cov}}
\newcommand{\diam}{\mathrm{diam}}
\makeatletter
\newsavebox\myboxA
\newsavebox\myboxB
\newlength\mylenA
\newcommand*\xoverline[2][0.70]{%
    \sbox{\myboxA}{$\m@th#2$}%
    \setbox\myboxB\null% Phantom box
    \ht\myboxB=\ht\myboxA%
    \dp\myboxB=\dp\myboxA%
    \wd\myboxB=#1\wd\myboxA% Scale phantom
    \sbox\myboxB{$\m@th\overline{\copy\myboxB}$}%  Overlined phantom
    \setlength\mylenA{\the\wd\myboxA}%   calc width diff
    \addtolength\mylenA{-\the\wd\myboxB}%
    \ifdim\wd\myboxB<\wd\myboxA%
       \rlap{\hskip 0.5\mylenA\usebox\myboxB}{\usebox\myboxA}%
    \else
        \hskip -0.5\mylenA\rlap{\usebox\myboxA}{\hskip 0.5\mylenA\usebox\myboxB}%
    \fi}
\makeatother

\newtheorem{theorem}{Theorem}[section]
\newtheorem{corollary}[theorem]{Corollary}
\newtheorem{conjecture}[theorem]{Conjecture}
\newtheorem{proposition}[theorem]{Proposition}
\newtheorem{lemma}[theorem]{Lemma}
\newtheorem{definition}[theorem]{Definition}
\newtheorem{example}[theorem]{Example}
\newtheorem{remark}[theorem]{Remark}
\newtheorem{case}{Case}[section]
\newtheorem{condition}{Condition}[section]
\newcommand{\proof}{\noindent {\it Proof:} }

\newcommand*\samethanks[1][\value{footnote}]{\footnotemark[#1]}

\title{{\bf\Large Normal approximation for fire incident simulation using permanental Cox processes}}
\author{Dawud Thongtha\thanks{Research partially supported by the TSRI Fundamental Fund 2020 (Grant Number: 64A306000047).} \ and Nathakhun Wiroonsri\samethanks \\ \normalsize Mathematics and Statistics with Applications Research Group, \\
	\normalsize Department of Mathematics,
	\normalsize King Mongkut's University of Technology Thonburi}

\maketitle

\date{}

\begin{abstract}   
	Estimating the number of natural disasters benefits the insurance industry in terms of risk management. However, the estimation process is complicated due to the fact that there are many factors affecting the number of such incidents. In this work, we propose a Normal approximation technique for associated point processes for estimating the number of natural disasters under the following two assumptions: 1) the incident counts in any two distinct areas are positively associated and 2) the association between these counts in two distinct areas decays exponentially with respect to distance outside some small local neighborhood. Under the stated assumptions, we extend previous results for the Normal approximation technique for associated point processes, i.e., the establishment of non-asymptotic $L^1$ bounds for the functionals of these processes  \cite{Wir19}. Then we apply this new result to permanental Cox processes that are known to be positively associated. Finally, we apply our Normal approximation results for permanental Cox processes to Thailand’s fire data from 2007 to 2020, which was collected by the Geo-Informatics and Space Technology Development Agency of Thailand.
\end{abstract}


\textbf{JEL}: C020 Mathematical Methods

\section{Introduction}
\indent Probability and statistical models have been widely used in most business sectors, including the insurance business. They are major tools for estimating loss, claim severity, claim count and claim probabilities, as these factors significantly affect their business operations. Since claim occurrences are used in determining policy premiums and are usually not predictable, claim simulation and prediction are some of the main tools for handling the situation. Simulation and prediction techniques for the insurance industry have attracted researchers’ attention for quite some time (see \cite{Jessen2011}, \cite{Godfred14} and \cite{Adrea18} and references therein). 

\indent The loss can be estimated as the product of the frequency of risk events and some overall measure of severity. Therefore, approximating the number of occurrences of risk events and claim probabilities is key to successful loss estimation. These related topics have been studied and developed in various directions. For instance, the works \cite{Lamb2005} and \cite{Caroline13} improved statistical and stochastic methods for estimating the probability and frequency of flood events. Machine learning techniques have also been applied to estimate claims and reserves (see \cite{Quan2018}, \cite{Wuthrich2018} and \cite{Bartl2020}, for example). 

\indent A point process (sometimes called a random point field) is one of the most powerful tools in probability. These processes, which have been widely developed from their root, namely, the Poisson point process, have been applied in several areas, including insurance (see\cite{Xu2011}, \cite{Gustav2015} and \cite{Albrecher2021}). 
Intuitively, a point process is an appropriate tool for modeling the occurrences of risk events and insurance claims, as it explains random natural phenomenon across space and time. Many researchers have proposed different approaches based on point processes for modeling natural disasters, especially fires. In 2007, in \cite{Schoenberg2007}, it was shown that a simple point process outperforms the Burning Index  in predicting wildfire incidents in Los Angeles County. In 2011, the work \cite{Xu2011} extended the model in \cite{Schoenberg2007} by considering relevant covariates, such as historical spatial burn patterns and wind direction. Also, in 2015, the work \cite{Gustav2015} used a spatial point process to study the risk associated with insurance customers using geographical information systems.

\indent Some point processes, such as the classic Poisson point process, have the property that the intensities in the relevant areas are independent. However, when considering natural disasters or insurance claims, a dependency structure must be introduced. In 2009, the work \cite{Keef2009} proposed a copula approach to finding the joint probability distribution for hydrological variables and then using it to study the spatial dependence in extreme river flows and precipitation. 
When exact distributions are unknown, approximate distributions can be used to estimate probabilities of risk events. For example, the work \cite{Caroline13} applied a theoretical result related to an approximation of a conditional distribution in \cite{Heffernan2007} to develop a new method for estimating the probability of widespread flood events. In our work, we use a Normal approximation for associated point processes to estimate the total number of fire incidents. When approximating distributions, the major issue is dependency. Stein's method (\cite{Stein72}) is one of the best tools for handling dependency. This well-known method for approximating limiting distributions has the additional benefit of providing non-asymptotic bounds for these distributions . In this work, we use the concept of local dependency, which can be handled by Stein's method, as shown in \cite{ChenShao2004}.

This work is divided into two parts. First, we extend the theoretical results in \cite{Wir19} to provide non-asymptotic $L^1$ bounds for the functionals of associated point processes under more general assumptions. These assumptions are motivated by the nature of natural disasters, as discussed later in this section. We also apply our main theoretical results to the permanental Cox process detailed in Section \ref{sec:def}, which is known to be positively associated. In the second part, we simulate Thailand’s fire incidents using permanental Cox processes based on the satellite data on Thailand ‘s fire occurrences collected by the Geo-Informatics and Space Technology Development Agency (GISTDA) from 2007 to 2020. Finally, we check our Normal approximation results from the first part using our simulation results in the second part.   

Stein's method has been used to obtain bounds for point processes under various assumptions and Poisson processes (see \cite{Bar88}, \cite{BB92}, \cite{CX04}, \cite{CX06} and \cite{CX11} for examples). However, association properties did not appear in the context of Normal approximations for point processes until the work \cite{PDL17} proved a central limit theorem for their functionals under certain conditions and the work \cite{Wir19} obtained rates of convergence using Stein's method under a stronger assumption. 

\indent Next, we state our assumptions regarding natural disasters. Based on the literature, natural disaster incidents in a particular area seem to affect such incidents in another area. For instance, according to \cite{Xu2011}, wind may cause fire to extend to nearby areas . Thus, we are inspired to claim that for some natural disasters, such as fires, there is some positive relation between the incident counts for two distinct areas. Moreover, we also claim that the further apart the two areas are, the weaker the relation should be. Therefore, we make the following two assumptions in our work:

\begin{enumerate}[label=(A\arabic*)] 
	\item The chances of incidents in any two areas are positively associated. Thus, incident counts in any two distinct areas are positively associated. 
	\item The association between the incident counts in two distinct areas decays exponentially with respect to the distance between the areas. However, the decay is assumed to begin outside some small local neighborhood.
\end{enumerate}

Though our assumptions are quite intuitive, to the best of our knowledge, we are unaware of any other works that use a dependent structure with respect to location to model natural disaster incidents.

The remainder of this work is organized as follows. We provide some necessary background and definitions in Section \ref{sec:def}. Then the theoretical results are stated and proved in Section \ref{sec:main}. Section \ref{sec:app} is devoted to the simulation of fires in Thailand using the dataset from GISTDA and checking our normal approximation results based on this simulation. Finally, a conclusion is provided in Section \ref{sec:sum}. 
\section{Some background and definitions} \label{sec:def}

We devote this section to a review of the theoretical background of point processes, permanental Cox processes and Stein's method.

\subsection{Point and permanental Cox processes}

A point process is a collection of random points on some mathematical space, such as Euclidean space $\mathbb{R}^{d}$. As mentioned earlier, point processes are popular tools in the insurance industry (see \cite{Albrecher2021}, \cite{Gustav2015} and \cite{Xu2011}).

Let $X \subset \mathbb{R}^d$ with $d \in \mathbb{N}$ be a point process. For $\text{card}(A)$, the cardinality of $A \subset \mathbb{R}^d$, let
\beas
N(A) = \text{card}(X \cap A).
\enas 

The point process $X$ is said to be \textit{simple} if $N(\{{\bf a}\}) \in \{0,1\}$ a.s. for all ${\bf a} \in \mathbb{R}^d$. Also, the process $X$ is said to be \textit{locally finite} if it takes values in 
\beas
\Omega = \left\{x \subset \mathbb{R}^d:\text{card}(x \cap B)<\infty,\forall \text{ bounded } B\subset \mathbb{R}^d\right\}.
\enas

The process $X$ is said to be \textit{negatively associated} if for all coordinate-wise increasing functions $\psi:\mathbb{N}^k\rightarrow \mathbb{R}$ and $\phi:\mathbb{N}^l\rightarrow \mathbb{R}$ and for all families of pairwise disjoint Borel sets $\{A_i \mid 1\le i \le k\}$ and $\{B_j \mid 1\le j \le l\}$  such that
\bea \label{disj}
(\displaystyle\cup_{i}A_i) \cap (\displaystyle\cup_{j}B_j) = \emptyset,
\ena
we have
\begin{multline*}
	\E\left[\psi\left(N(A_1),\ldots,N(A_k)\right)\phi\left(N(B_1),\ldots,N(B_l)\right)\right] \\
	\le \E\left[\psi\left(N(A_1),\ldots,N(A_k)\right)\right] \E \left[\phi\left(N(B_1),\ldots,N(B_l)\right)\right] .
\end{multline*}
Similarly, the point process is said to be \textit{positively associated} if the above inequality is reversed and the families of Borel sets are not necessarily assumed to satisfy \eqref{disj}. In addition, a point process is said to be \textit{associated} if it is either negatively or positively associated. 

Next, we state the definition of the nth order intensity functions of point processes with respect to Lebesgue measure. 
Let $n \in \mathbb{N}$ and $X \in \Omega$. If there exists a non-negative function $\rho_n : \left(\mathbb{R}^d\right)^n \rightarrow \mathbb{R}$ such that
\beas
\E\left[\sum_{\substack{{\bf x}_1,\ldots,{\bf x}_n \in X \\ \text{all distinct}}}f({\bf x}_1,\ldots,{\bf x}_n)\right] 
= \int_{\left(\mathbb{R}^d\right)^n} f({\bf x}_1,\ldots,{\bf x}_n) \rho_n({\bf x}_1,\ldots,{\bf x}_n) d{\bf x}_1 \ldots d{\bf x}_n
\enas
for all locally integrable functions $f : \left(\mathbb{R}^d\right)^n \rightarrow \mathbb{R}$, then $\rho_n$ is called the \textit{nth order intensity function} of $X$.
Now, for ${\bf x},{\bf y} \in \mathbb{R}^d$, let
\bea \label{ddef}
D({\bf x},{\bf y}) = \rho_2({\bf x},{\bf y})-\rho_1({\bf x})\rho_1({\bf y}).
\ena
It follows that
\beas
\Cov(N(A),N(B)) = \int_{A \times B} D({\bf x},{\bf y}) d{\bf x} d{\bf y}.
\enas

Cox processes are well-known point processes. They are considered to be generalizations of the Poisson point processes for which the intensity is a random measure.   
A permanental Cox process \cite{MM05} is the Cox process with intensity functions
\bea \label{coxint}
\rho_n({\bf x}_1,\ldots,{\bf x}_n) = \E \prod_{i=1}^n \Lambda({\bf x}_i),
\ena
where $\Lambda $ is a random measure defined as
\bea \label{coxgauss}
\Lambda({\bf x}_i) = Y_1^2({\bf x}_i)+\cdots+Y_l^2({\bf x}_i),
\ena 
where the $Y_1,\ldots,Y_l$ are $l$  independent, zero-mean, real-valued Gaussian random fields with covariance function $C$.

In this work, we consider the case where the Gaussian random fields are stationary; hence, $C$ only depends on $r = |{\bf x}_i-{\bf x}_j|_\infty$, where $|\cdot |_\infty$ denotes the vector max norm. Also, we consider just the specific case when 
\bea \label{covpat}
C(r) = \kappa_c e^{-\lambda_c r},
\ena
where $\kappa_c$ and $\lambda_c$ are constants. This process has been used, for instance, in \cite{MY06}, \cite{YMM12} and \cite{MESIT18}. It is known to be positively associated (see \cite{Eis14} and \cite{LSY20}).

\subsection{Stein's bound for local dependent random variables}
Stein's method, introduced by Charles Stein \cite{Stein72} in 1972, is a widely known technique for finding non-asymptotic bounds for approximations of probability distributions. It was motivated by the idea that $W$ has the standard Normal distribution
if and only if 
\beas
\E W f(W) = \E f'(W)
\enas
for all absolutely continuous functions $f$ with $\E |f'(W)| < \infty$. This identity leads to the differential equation
\bea \label{steineq}
h(w) - \E h(Z) = f'_h(w)-wf_h(w),
\ena
where $Z$ is a standard Normal random variable and $h$ is a test function. If $h \in \mathcal{H}$ for some Borel set $\mathcal{H}$ and we replace $w$ by a random variable $W$, the error in the distributional approximation of $W$ by $Z$ on $\mathcal{H}$ can be bounded by obtaining the non-asymptotic bound of the expectation of the supremum of the right-hand side of \eqref{steineq} over $\mathcal{H}$. In general, doing so is much simpler than computing the left-hand side of \eqref{steineq} directly. Taking $h \in \mathcal{H}_{1}$, where $\mathcal{H}_{1}=\{h \mid |h(x)-h(y)|\leq |x-y|\}$, we obtain
\beas 
d_{1} (\mathcal{L}(W),\mathcal{L}(Z)) = \sup_{h \in \mathcal{H}_{1}} |\E h(W)-\E h(Z)|.
\enas 
This distance is known as $L^1$ or the Wasserstein distance. Stein's method has been used in various applications, and it is one of the best ways to handle dependent situations (see \cite{CGS11} and \cite{Ross11}). One of the classic dependent cases, handled by Stein's method, is the local dependent structure introduced in \cite{ChenShao2004}. A collection of random variable $X_1,\ldots,X_n$ has dependency neighborhoods $N_{i} \subseteq \{1,2,\ldots, n\}, i=1,2,\ldots, n $ if $X_{i}$ is independent of $X_{j}$ for all $j \notin N_i - \{i\}$. Next, we state a version of the local dependence bound that appeared in the note \cite{Ross11}.

\begin{theorem}[\cite{Ross11}] \label{local}
	Let $X_1,\ldots,X_n$ be random variables such that $\E[X_i^4] < \infty$, $\E X_i = 0$, $\sigma^2 = \Var(\sum_i X_i)$, and define $W = \sum_i X_i/\sigma$. Let the collection $(X_1,\ldots,X_n)$ have dependency neighborhoods $N_i$, $i=1,\ldots,n$, and also define $D = \max_{1\le i \le n}|N_i|$. Then for $Z$, a standard Normal random variable, we have
	\beas
	d_1\big({\cal L}(W),{\cal L}(Z)\big) \le \frac{D^2}{\sigma^3}\sum_{i=1}^n\E|X_i|^3 + \frac{\sqrt{28}D^{3/2}}{\sqrt{\pi}\sigma^2}\sqrt{\sum_{i=1}^n\E[X_i^4]}.
	\enas 
\end{theorem}
We apply Theorem \ref{local} to obtain our main results in the next section.

\section{Main results}  \label{sec:main}
In this section, we state and prove our main result, which is extended from \cite{Wir19} and then apply it to permanental Cox Processes. In this work, we relax the assumption in \cite{Wir19} that the covariance of the second-order intensity decays exponentially and assume that it decays exponentially outside of some local neighborhood. 
Let $X \in \Omega$ be an associated point process and
\bea \label{Ydef}
Y_{\bf i} = f_{\bf i}(X \cap C_{\bf i}) - \E f_{\bf i}(X \cap C_{\bf i}), {\bf i} \in \mathbb{Z}^d,
\ena
where the $f_{\bf i}:\Omega \rightarrow \mathbb{R}$ are real-valued measurable functions and the $C_{\bf i}$ are defined as \mbox{$d$-dimensional} unit cubes centered at ${\bf i}$. Note that the union of $C_{\bf i}$ forms a covering of $\mathbb{R}^d$. 

We let ${\bf 1} \in \mathbb{Z}^d$ denote the vector with all components $1$, and write inequalities such as ${\bf a} < {\bf b}$ for vectors ${\bf a} , {\bf b} \in \mathbb{R}^d$ when they hold component-wise. In this work, we consider 
\bea \label{blockdef}
S_{{\bf k}}^n = \sum_{{\bf i} \in B_{{\bf k}}^n} Y_{{\bf i}} \qmq{where}
B_{{\bf k}}^n = \left\{ {\bf i} \in \mathbb{Z}^d: {\bf k} \le {\bf i} < {\bf k}+n{\bf 1} \right\}.
\ena   

%%%%%%%%%%%%%%%%%%%%%%%%%%%%%%%%%%%%%%%%%%%%%%%%%%%%%%%%%%%%%%%%%%%%%%%%%%%5
The work \cite{PDL17} obtained a CLT for the sum above with $B_{{\bf k}}^n$ replaced by any sequence of strictly increasing finite domains of $\mathbb{Z}^d$. Let $C_{\bf i}$ be any $d$-dimensional cube centered at $x_{\bf i} = R \cdot {\bf i}$ with fixed $R>0$ and with fixed side length $s \ge R$. We state our main result for $R=s=1$. We then add a remark after the theorem to discuss how to generalize this result. In the following, for $k \in \mathbb{N}$, we denote $\left\|X\right\|_{k}= \left\lbrace \mathbb{E}|X|^{k}\right\rbrace^{\frac{1}{k}}$.
%%%%%%%%%%%%%%%%%%%%%%%%%%%%%%%%%%%%%%%%%%%%%%%%%%%%%%%%%%%%%%%%%%%%%%%%555

\begin{theorem} \label{localplus}
	For $d \in \mathbb{N}$, let $X$ be a locally finite simple associated point process on $\mathbb{R}^d$. For ${\bf k} \in \mathbb{Z}^d$, let $S_{{\bf k}}^n$ be as in \eqref{blockdef}, with $Y_{\bf i}$ given in \eqref{Ydef}, $R=s=1$ and $\sigma_{n,{\bf k}}^2 = \Var(S_{{\bf k}}^n)$. Assume that the following conditions are satisfied:
	\begin{enumerate}[label=(\alph*)]
		\item The first two intensity functions of $X$ are well defined; \label{a}
		\item $\sup_{{\bf i} \in \mathbb{Z}^d} \left\|Y_{\bf i}\right\|_{4} = M < \infty$; \label{b}
		\item $\sup_{|{\bf x}-{\bf y}|_\infty \ge r} D({\bf x},{\bf y}) \le \kappa e^{-\lambda r}$ for some $\kappa,\lambda > 0$ and $r>Kn^{\frac{2}{3}\left(\frac{4d-1-1/d}{4d+2} \right)} $ with $K>0$; \label{c}
		\item $\sigma_{n,{\bf k}}^2 \ge \gamma n^d$ for some $\gamma>0$. \label{d}
	\end{enumerate}
	Then, for the standard Normal random variable $Z$,
	\beas 
	d_1\left(\mathcal{L}\left(\frac{S_{{\bf k}}^n}{\sigma_{n,{\bf k}}}\right),\mathcal{L}(Z)\right) &\le&
	\left(\frac{\sqrt{28}M^2(2K)^{3d/2}}{\gamma \sqrt{\pi}}+ C_{1,d,M,\kappa,\gamma}\right)\frac{1}{n^{d/(4d+2)}}+ \frac{M^3(2K)^{2d}}{\gamma^{3/2}n^{d/6-1/(6d+3)}} \nn\\
	&& \hspace{5pt}+ \frac{C_{2,d,M,\kappa,\gamma} n^{d(4d+1)/(6d+3)}}{\exp\left(\theta_{d,M,\kappa,\gamma}n^{d/(4d+2)}\right)}+ \frac{C_{3,d,M,\kappa,\gamma} n^{7d/6}}{\exp\left(2\theta_{d,M,\kappa,\gamma}n^{d/(4d+2)}\right)},
	\enas
	where $d_1$ is the $L^1$ distance,\\
	\beas
	\theta_{d,M,\kappa,\gamma} =\frac{\lambda}{3}\left( \frac{\sqrt{2\gamma}\kappa^{1/3}\left((4\mu_{\lambda}+2\nu_{\lambda})^d-\left(2\nu_{\lambda}\right)^d\right)}{18^{d+1}\sqrt{\pi}dM}\right)^{1/(2d+1)},
	\enas
	\beas
	&&C_{1,d,M,\kappa,\gamma} = \left(\frac{9\cdot 36^d M^{4d+3}\left((4\mu_{\lambda}+2\nu_{\lambda})^d-\left(2\nu_{\lambda}\right)^d\right)^{2d} }{\gamma^{2d+(3/2)}\pi^d}\right)^{1/(2d+1)} \\ 
	&& \hspace{100pt} \times \left(\frac{1}{(2d)^{2d/(2d+1)}}+2(2d)^{1/(2d+1)}\right) ,
	\enas
	\beas 
	C_{2,d,M,\kappa,\gamma} = \frac{3\cdot 6^d \kappa^{1/3} M^2 \theta_{d,M,\kappa,\gamma}^{4d/3}}{\sqrt{\pi}\gamma}, \ \ \ 
	C_{3,d,M,\kappa,\gamma} = \frac{2^{d+1} \kappa^{2/3} M}{\sqrt{\gamma}},
	\enas
	and
	\beas 
	\mu_\lambda = \frac{e^{\frac{2\lambda}{r}}}{\left(e^{\frac{\lambda}{r}}-1\right)^2}, \ \ \ \nu_\lambda = \frac{e^{\lambda/r}}{\left(e^{\frac{\lambda}{r}}-1\right)^2}.
	\enas
\end{theorem}

\proof
We prove this theorem by applying Theorem \ref{local} stated above and Theorem 3.1 in \cite{Wir19}, handling local dependence and non-local dependence separately. First, we address local dependence only by assuming that
	 $\Cov(Y_{\bf i},Y_{\bf j})=0$ for $|{\bf i}-{\bf j}|_\infty> Kn^{\frac{2}{3}\left(\frac{4d-1-1/d}{4d+2} \right)}$.
Invoking Theorem \ref{local} with 
	$D = (2K)^{d} n^{\frac{2}{3}\left(\frac{4d^{2}-d-1}{4d+2} \right)}$, 
 Assumption \ref{b} that $\E|Y_{\bf j}|^3\le M^3$ and \mbox{$ \E[Y_{\bf j}^4] \le M^4$} and using Assumption \ref{d},
we have that the $L^1$ distance is bounded by
\bea \label{bdd1}
\frac{\sqrt{28}M^2 (2K)^{3d/2}}{\gamma \sqrt{\pi}n^{d/(4d+2)}} + \frac{M^3 (2K)^{2d}}{\gamma^{3/2}n^{d/6-1/(6d+3)}}.
\ena

We now assume that Assumption \ref{c} is true for all $r>0$. Invoking Theorem 3.1 in \cite{Wir19}, we have that the $L^1$ distance is bounded by 
\bea \label{bdd2}
\frac{C_{1,d,M,\kappa,\gamma}}{n^{d/(4d+2)}} + \frac{C_{2,d,M,\kappa,\gamma} n^{d(4d+1)/(6d+3)}}{\exp\left(\theta_{d,M,\kappa,\gamma}n^{d/(4d+2)}\right)} + \frac{C_{3,d,M,\kappa,\gamma} n^{7d/6}}{\exp\left(2\theta_{d,M,\kappa,\gamma}n^{d/(4d+2)}\right)}
\ena

The distance $d_1\big({\cal L}(W_{{\bf k}}^n),{\cal L}(Z)\big)$ is bounded by the sum of \eqref{bdd1} and \eqref{bdd2} because the local neighborhood covariance terms and the remaining covariance terms have been handled by \eqref{bdd1} and \eqref{bdd2}, respectively.
\bbox

\begin{remark}  
	\begin{enumerate}
		\item The associated assumption and Assumption \ref{c} are motivated by (A1) and (A2) in the introduction, respectively.
		\item The size of the local neighborhood is $O\left(n^{d/(4d+2)}\right)$, which is flexible and can be increasing in $n$. 
		\item The above theorem can be extended to the case where $R$ and $S$ are greater than one by following the same proof. The bound will end up with a larger constant. 
	\end{enumerate}
\end{remark}
Note that we added the assumption that \ref{c} is true for all $r>0$ in the above proof after the local covariance terms had been handled. Although doing so may make the constant larger, the rate of convergence is unaffected. 

Next, we apply the above theorem to the permanental Cox process on $\mathbb{R}^d$, which is known to be positively associated. In this theorem, we consider a function $f: \Omega \rightarrow \mathbb{R}$ defined by
\bea \label{fdef}
f(Y) = \sum_{S \subset Y} g(S)\mathbf{1}_{|S|=p},
\ena
where $g$ is a bounded function supported on sets $S$ having exactly $p$ elements such that $g(S)=0$ when $diam(S)> \tau$ for some fixed $\tau>0$ and $p \in \mathbb{N}$. Here, we denote $diam(S)=\sup_{x,y \in S} \left|x-y\right|_\infty$. Also, we focus specifically on the case that $p=1$. We note here that if $p=1$ and $g(S) = 1$, $f(X \cap \Lambda_n)$ is the number of points in $X \cap \Lambda_n$ or $N(\Lambda_n)$, where $\Lambda_n \subset \mathbb{R}^d$. The result is stated as follows:

\begin{theorem} \label{mainthm2}
	Let $n,k \in \mathbb{N}$ and $X$ be a permanental Cox process with intensity functions as in \eqref{coxint} and $\Lambda({\bf x}_i)$ as in \eqref{coxgauss}, where the $Y_i$, $i=1,\cdots,k$ are independent, mean-zero, real-valued Gaussian random fields. Let $f$ be defined as in \eqref{fdef} with $g$ be bounded and $g(S)=0$ when $diam(S)> \tau$ for some fixed $\tau>0$ and $p=1$. Letting $\sigma_n^2 = \Var(f(X \cap \Lambda_n))$ with $\Lambda_n \subset \mathbb{R}^d$ be such that $|\Lambda_n|=n^d$, assume that 
	\bea \label{as1}
	\sup_{|{\bf x}-{\bf y}|_\infty \ge r} \Cov \left(\Lambda({\bf x}),\Lambda({\bf y})\right) \le O(e^{-\lambda r}),   
	\ena
	for some $\lambda >0$ and $r>Kn^{d/(4d+2)}$ with $K>0$, and
	\bea \label{as2}
	\liminf_n \sigma_n^2/n^d >0.
	\ena
	Then for $W_n = (f(X \cap \Lambda_n)-\E f(X \cap \Lambda_n))/\sigma_n$, there exists $C>0$ such that
	\bea \label{mainbound2}
	d_1\left(\mathcal{L}(W_n),\mathcal{L}(Z)\right) \le \frac{C}{n^{d/(4d+2)}}, 
	\ena
	where $Z$ is a standard Normal random variable. 
\end{theorem}

\proof To prove the theorem, we follow the same argument as that used in the proof of Theorem 4.3 in \cite{Wir19}. It requires the use of Theorem \ref{localplus}, where the assumptions \eqref{as1} and \eqref{as2} are needed. Because Theorem 4.3 of \cite{Wir19} is for determinantal point processes, the variance condition used in its proof cannot be used here; thus, the proofs differ at this point. Therefore, it is sufficient to show that
\beas
\Var\left(\sum_{S \subset X \cap \Lambda_n} g(S)\mathbf{1}_{|S|=p}\right) = O(\Lambda_n),
\enas
which can be proved by adapting the proof of Lemma B.6 in \cite{PDL17} for the case where $X$ is a permanental Cox process. It is sufficient to verify that the term (B.4) in that proof is bounded for $k=0,\ldots,p$. For the term $k=0$, we show that the intensity $\rho_{2p}$ is bounded as follows: 
\beas
|\rho_{2p}({\bf x}_1,\cdots,{\bf x}_{2p})| &=& \E \prod_{j=1}^{2p} \Lambda({\bf x}_j) =\E \prod_{j=1}^{2p} \left(Y^2_1({\bf x}_j)+\cdots+Y^2_l({\bf x}_j)\right) \\
&=& \sum_{1\le i_1,\ldots,i_{2p} \le l}\E \left[Y^2_{i_1}({\bf x}_1)Y^2_{i_2}({\bf x}_2)\cdots Y^2_{i_{2p}}({\bf x}_{2p})\right] \\
&\le& \sum_{1\le i_1,\ldots,i_{2p} \le l} \left(\E \left[Y^{2p}_{i_1}({\bf x}_1)\right] \E \left[Y^{2p}_{i_2}({\bf x}_2)\right] \cdots \E \left[Y^{2p}_{i_{2p}}({\bf x}_{2p})\right] \right)^{1/p} \\
&\le& l^{2p} \left(M(2p-1)!!\right)^2,
\enas
where $M=\max_{1\le i\le l,{\bf x} \in \mathbb{R}^d}\sigma^{2}_i({\bf x})$, with $\sigma^2_i({\bf x})=\Var(Y_i({\bf x}))$. Note that we have used H\"{o}lder's inequality and the fact that $Y_i$ is Normal and its central moment is\begin{center}
	 $\E |Y_{i}({\bf x})|^{2p}= \sigma^{2p}_i({\bf x})(2p-1)!!$.
\end{center}  As we only consider the case $p=1$, the term $k= 1$ is obvious. Thus, the order of the variance is verified. 
\bbox

\begin{remark}  
	We remark here that in the proof of Theorem \ref{mainthm2}, we set $M = \max_{1\le i\le l,{\bf x} \in \mathbb{R}^d}\sigma_i^{2}({\bf x})$, which could be large. However, in our application, we later set the $\sigma_i^2({\bf x})$ to be equal for all $i$ and $\sigma_i^2({\bf x}) = m({\bf x})/l$, where $m({\bf x})$ is the maximum number of fires in area ${\bf x}$, as derived from the historical data, which is assumed to be bounded.  
\end{remark}

\section{Application to the Thai fire dataset} \label{sec:app}

In this section, we use our main results to simulate Thailand’s fires via permanental Cox processes. We use the GISTDA Thai fire dataset, collected from 2007 to 2020 by satellite. The dataset consists of the latitudes and longitudes of all fire incident locations. Recall that we claim that the fire counts from two distinct areas should be positively correlated. Moreover, we claim that their covariance begins to decay exponentially outside some small local neighborhood. We split this section into three subsections. We first explore the dataset and check to see that Assumptions (A1) and (A2) are not contradicted by the data in the first subsection. The second subsection is devoted to applying our main results to fire simulations. In the last subsection, we evaluate the effect of varying the decay parameter for the permanental Cox processes.

\subsection{Exploring the dataset} \label{sec:explore}

\begin{figure}[H]
	\begin{center}
		\includegraphics[width=3.5in]{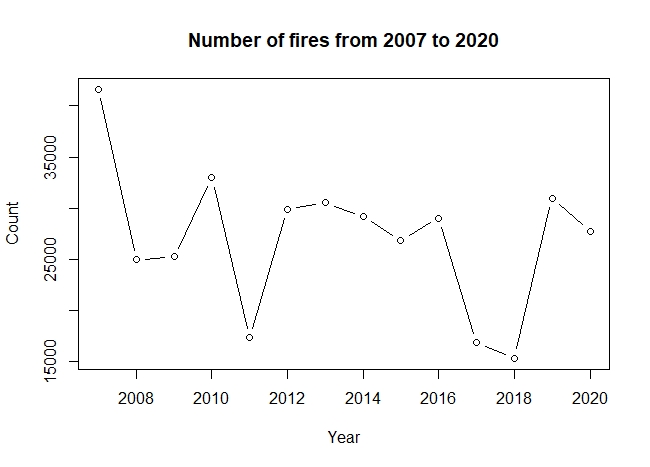}
	\end{center}
	\caption{The total number of yearly  fires for the period 2007-2020}
	\label{total}      
\end{figure}

\begin{figure}[H]
	\begin{center}
		\includegraphics[width=5in]{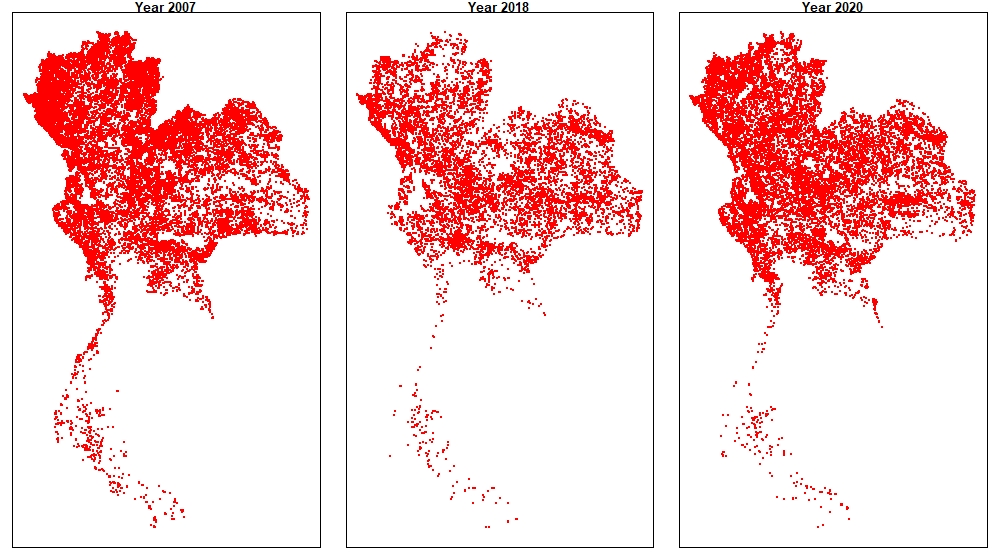}
	\end{center}
	\caption{The actual fire locations in 2007, 2014 and 2020}
	\label{fig:realplot}      
\end{figure}

In this subsection, we first explore the dataset and verify that it does not contradict the assumptions in Theorem \ref{mainthm2}. We follow the process outlined in Figure \ref{fc41}. Figure \ref{total} shows the total number of fires each year in Thailand, where the mean and the standard deviation were  27,048 and 6,980.94, respectively. Figure \ref{fig:realplot} shows the locations of fires in 2007, 2014 and 2020, respectively, from left to right.

\begin{figure}[H]
	\centering
	\scalebox{0.8}{
		\begin{tikzpicture}[node distance=1.8cm,
			every node/.style={fill=white, font=\sffamily}, align=center]
			\node (Data)            		[base, fill=red!30]	{Fires data 2007 - 2020};
			\node (Assumption)      		[base, below of=Data, fill=green!40]	{Verifying the assumptions \\ (Theorem \ref{mainthm2})};
			\node (Descriptive)     		[base, right of=Data, xshift=4cm, fill=blue!20]	{Descriptive statistics};
			\node (Condition22)			[base, below of=Assumption, xshift=-4cm, node distance=2cm, fill=green!25]		{Assumption \eqref{as1}};
			\node (Condition23) 		[base, below of=Assumption, xshift=4cm, node distance=2cm, fill=green!25]		{Assumption \eqref{as2} \\ The variance order};
			\node (LocalNonDecay)		[base, below of=Condition22, xshift=-2.5cm, node distance=2cm, fill=green!10]	{Local neighborhood};
			\node (ExponentialDecay)	[base, below of=Condition22, xshift=2.5cm, node distance=2cm, fill=green!10]	{Exponential decay};

			\draw[->]		(Data) -- (Descriptive);
			\draw[->]		(Data) -- (Assumption);
			\draw[->]		(Assumption) -- (0,-2.8) -| (Condition22);
			\draw[->]		(Assumption) -- (0,-2.8) -| (Condition23);
			\draw[->]		(Condition22) -- (-4,-4.8) -| (LocalNonDecay);
			\draw[->]		(Condition22) -- (-4,-4.8) -| (ExponentialDecay);
			\label{fc41}
		\end{tikzpicture}
	}
	\caption{Exploring the GISTDA fire dataset} \label{fc41}
\end{figure}
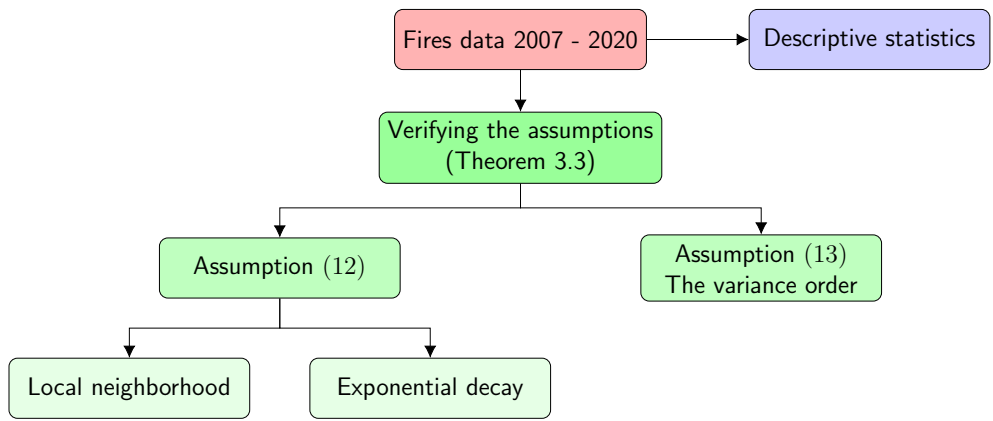

We now examine whether the fire dataset follows assumptions \eqref{as1} and \eqref{as2}. As the actual intensity function is unknown and cannot be obtained from the dataset, we can only check the covariances between fire counts in two distinct areas separated by a distance $r$. We first check to see if the covariance does not decay in a small local neighborhood. Then we further check to see if the covariance decays exponentially as the distance between areas increases. We also show the variances of the fire counts corresponding to $n$, with $n$ going from 10 to 645, where $645$ is the greatest distance possible in Thailand (1 unit equals to 0.01 degrees of latitude or longitude). Note that the total area of Thailand is 513,120 square kilometers, which equals $716^2$. Since 1 degree is about 111 kilometers, the total area is approximately $6.45^2$ square degrees, which is equivalent to $645^2$ square units when 1 unit is 0.01 degrees. Therefore, the largest $n$ possible in Thailand is 645.  

\begin{figure}[h]
	\begin{center}
		\includegraphics[width=5in]{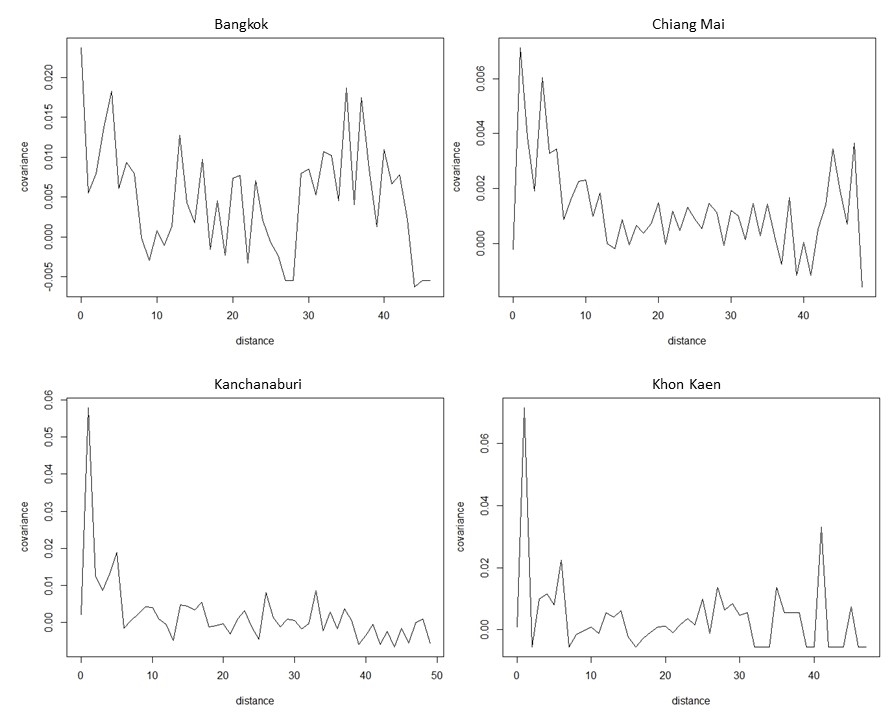}
	\end{center}
	\caption{Covariance decay patterns for Bangkok, Chaing Mai, Kanchanaburi and Khon Kaen with a unit distance of 0.002 degrees}
	\label{fig:bkk}      
\end{figure}

\begin{figure}[H]
	\begin{center}
		\includegraphics[width=3.5in]{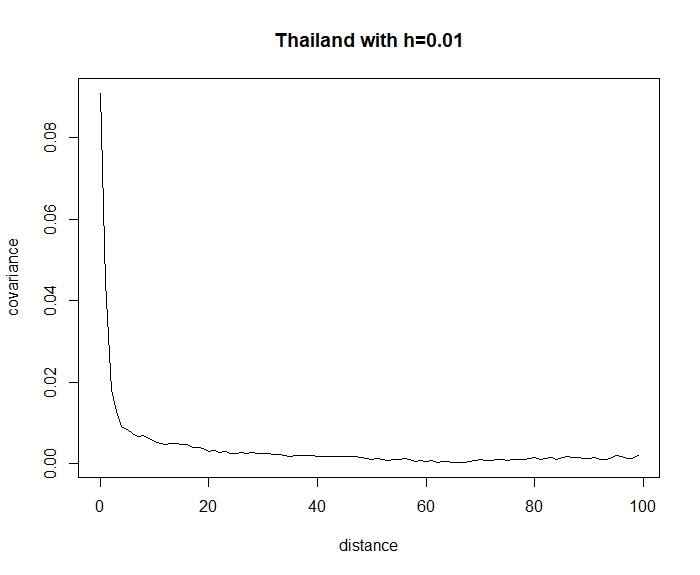}
	\end{center}
	\caption{Covariance decay pattern for Thailand with a unit distance of 0.01 degrees}
	\label{fig:thailand}      
\end{figure}

\begin{figure}[H]
	\begin{center}
		\includegraphics[width=3.5in]{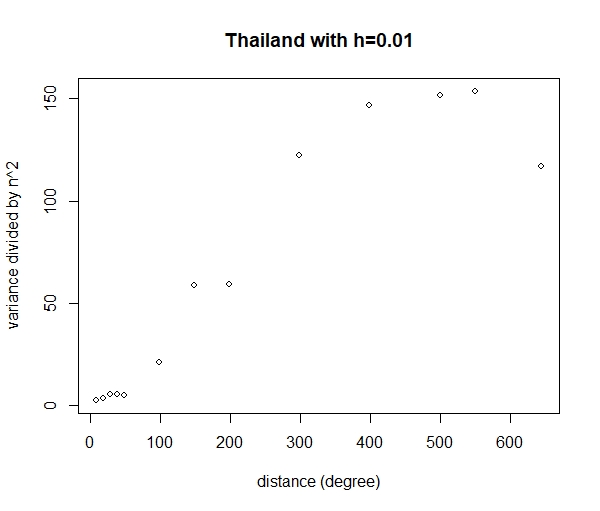}
	\end{center}
	\caption{Variance of the total number of fires divided by $n^2$ }
	\label{realvarrate}      
\end{figure}

As the entire country is too large to be considered a small local neighborhood, we specifically check four provinces from four regions, including Bangkok, Chaing Mai, Kanchanaburi and Khon Kaen. Note that we did not select a province from Southern Thailand, as the number of fires in this region was too low. Figure \ref{fig:bkk} shows no sign of exponential decay in the covariance pattern for these provinces when the unit distance is 0.002 degrees ($\sim 222$ meters). Figure \ref{fig:thailand} shows the covariance decay pattern based on distance for all of Thailand. We obtain an exponential decay rate of 0.15 when we use a unit distance of 0.01 degrees ($\sim 1.11$ kilometers). Figure \ref{realvarrate} shows the values of $\sigma_n^2/n^2$ for $n$ from 10 to 645. Obviously, assumption \ref{as2} is not contradicted.

\subsection{Simulation of Thai fires using the permanental Cox process} 

\begin{figure}[H]
	\centering
	\scalebox{0.8}{
		\begin{tikzpicture}[node distance=1.8cm,
			every node/.style={fill=white, font=\sffamily}, align=center]
			\node (Data)    [base, fill=red!30]	{Fires data 2007 - 2020};
			\node (EstPar1)	[base,below of=Data, fill=brown!40]	{Estimating permanental Cox process \\ parameters $\sigma^2, l$ in \eqref{coxgauss}};
			\node (EstPar2)	[base, below of=EstPar1, fill=brown!25]	{Estimating covariance decaying \\ parameters $\lambda_c$ and $K_c$ in \eqref{covpat}};
			\node (SimFires)	[base, below of=EstPar2, fill=brown!15]	{Simulating Fires corresponding to \\ the estimated parameters};
			\node (NormTest)	[base, left of=SimFires, xshift=-4.5cm, fill=blue!20]	{Testing normality of \\ the simulated fires};
			\node (VerAssump) [base, below of=SimFires, fill=green!30]	{Verifying assumptions \eqref{as1} and \eqref{as2}  \\ for the simulated fires};
			\node (ComputeL1)	[base, right of=SimFires, xshift=6cm, fill=orange!30] {Computing the $L_1$ distance \\ between the number of simulated fires \\ and a standard normal};
			\node (Compare)	[base, below of=ComputeL1, fill=yellow!20]	{Comparing the rate of convergence \\ to the bound in \eqref{mainbound2}};
			
			\draw[->]	(Data) -- (EstPar1);
			\draw[->]	(EstPar1) -- (EstPar2);
			\draw[->]	(EstPar2) -- (SimFires);
			\draw[->]	(SimFires) -- (NormTest);
			\draw[->]	(SimFires) -- (VerAssump);
			\draw[->]	(SimFires) -- (ComputeL1);
			\draw[->]	(VerAssump) -- (Compare);
			\draw[->]	(ComputeL1) -- (Compare);
			
		\end{tikzpicture}
	}
	\caption{Simulation of Thailand’s fires} \label{fc42}
\end{figure}
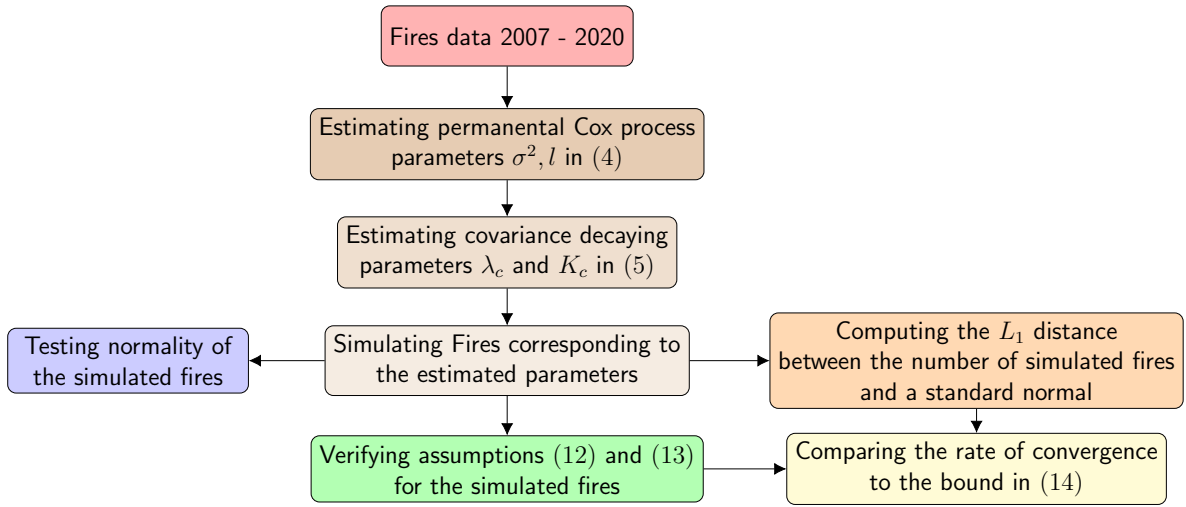

We devote this subsection to simulating fires in Thailand via the process outlined in Figure \ref{fc42}.

We first model Thailand’s fires using permanental Cox processes. Note that when we refer to the area ${\bf x}$, we are referring to a unit square that has ${\bf x}$ at the top left corner, where ${\bf x} \subset \mathbb{R}^2$ and 1 unit equals 0.01 degrees of latitude or longitude. We estimate the variances of $Y_j({\bf x}_{i})$ and $l$ from 
\beas
\Lambda({\bf x}_i) = Y_1^2({\bf x}_i)+\cdots+Y_l^2({\bf x}_i)
\enas
by using the first two sample moments. We set the means of the data and process to be equal and seek the closest variances. Writing $\sigma_j^2({\bf x}_i) = \Var(Y_j({\bf x}_i))$ for each ${\bf x}_i$ and using the fact that the $Y_j({\bf x}_{i})$ for $j=1,2,\ldots,l$ are zero-mean, independent Gaussian random variables, we set $\sum_{j=1}^l\sigma_j^2({\bf x}_i)$ equal to the sample mean of area ${\bf x}_i$ and $2\sum_{j=1}^l\sigma_j^4({\bf x}_i)$ equal to the sample variance of area ${\bf x}_i$. We then set $\sigma_j^2({\bf x}_i) = \sigma^2({\bf x}_i)$ for all $j$ for simplicity and solve for the closest $\sigma^2({\bf x}_i)$ and $l$. 

Denoting $m({\bf x}_i)$ and $v({\bf x}_i)$ as the mean and variance of the total fires in the unit cube at position ${\bf x}_i$ from the dataset, we have
\beas
\sigma^2({\bf x}_i)  = \frac{m({\bf x}_i)}{l} \text{ \ for all \ } i.
\enas
Then setting 
\beas
2l\sigma^4({\bf x}_i) = \sum_{j=1}^l\Var(Y^2_j({\bf x}_i)) = v({\bf x}_i),
\enas
and substituting $\sigma^2({\bf x}_i )$, we have
\beas
l  = \frac{2m^2({\bf x}_i)}{v({\bf x}_i)}.
\enas
For simplicity, we set $l = \max\left(1,\left\lfloor \frac{2m^2({\bf x}_i)}{v({\bf x}_i)}\right\rfloor\right)$. Although the $l$’s differ for different ${\bf x}_i$, we can choose the largest $l$. For a ${\bf x}_i$ with a smaller $l$, we set $\sigma^2({\bf x}_i)$ equal to zero for all the remaining terms.

To estimate the covariance function of the Gaussian processes from \eqref{covpat}, we need estimates of $\kappa_c$ and $\lambda_c$. Plugging in $r=0$, we have $\kappa_c = \sigma^2({\bf x}_i)$; thus, it suffices to estimate $\lambda_c$. We seek a $\lambda_c$ from a fine grid from 0 to 0.5 that results in the covariance decay rate closest to the 0.15 from the real dataset. Note that the decay rate of 0.15 is computed from the number of fires, whereas $\lambda_c$ is the decay rate of the covariance function for the Gaussian process. Therefore, it is not possible to set $\lambda_c = 0.15$ directly. 

\begin{figure}[H]
	\begin{center}
		\includegraphics[width=4in]{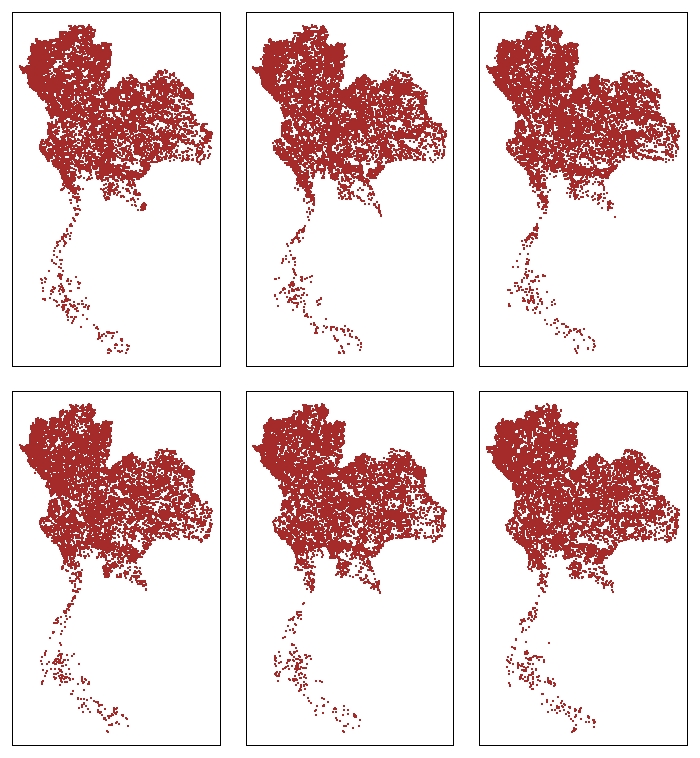}
	\end{center}
	\caption{This figure shows fire count simulation results using the estimated parameters and six iterations. Here, the simulation’s mean and standard deviation are 26,873 and 1078.15, respectively.}
	\label{fig:firesim}      
\end{figure}

By following the procedure above, we find the largest $l$, i.e., $l=8$, and closest $\lambda_c (0.1$). Figure \ref{fig:firesim} shows an example of a simulation with parameters estimated from the procedure above and 6 iterations.

\begin{figure}[h]
	\begin{center}
		\includegraphics[width=3.5in]{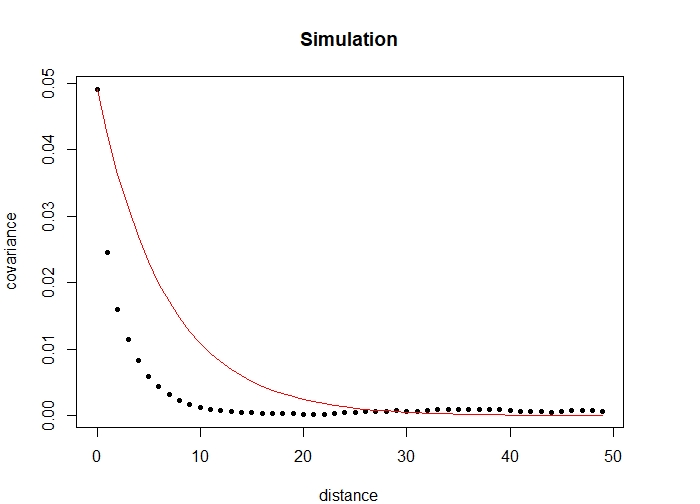}
	\end{center}
	\caption{The black dots show the covariances between the simulated fire counts for two distinct unit areas according to the distance between the areas. The red line shows exponential decay at a rate of 0.15.}
	\label{fig:simcovexp}      
\end{figure}

Next, we simulate fire counts for the whole country for 20 iterations and check that the covariance decays exponentially with distance. Since the number of areas in the entire country at a unit distance of 0.01 degrees from another area is extremely large, 20 is a reasonable number of iterations for the simulation. Figure \ref{fig:simcovexp} shows the covariance plot and the decay rate of 0.15 from the real dataset.

\begin{figure}[H]
	\begin{center}
		\includegraphics[width=3in]{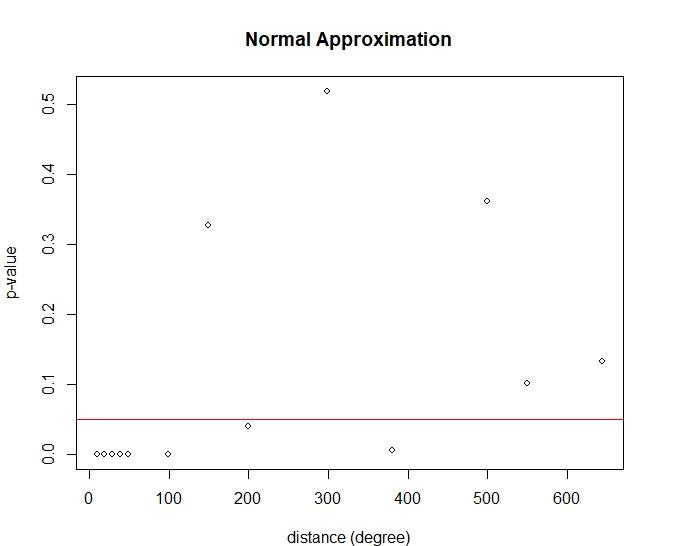}
	\end{center}
	\caption{The Shapiro–Wilk p-values for the fire counts}
	\label{fig:shapiro}      
\end{figure}

Then we simulate fires across the whole country for 500 iterations to check that the variance is on the order of $n^2$ and to compare the results with the Normal approximation from Theorem \ref{mainthm2}. As in Subsection \ref{sec:explore}, we take $n$ from 10 to 645, where $645$ units is the greatest distance possible in Thailand (1 unit equals 0.01 degrees of latitude or longitude). We preliminarily check the normality of the simulated fire counts using the Shapiro–Wilk test. Figure \ref{fig:shapiro} shows the Shapiro–Wilk p-values based on $n$, while the red horizontal line indicates the 0.05 significance level. We can see that when $n$ is large, the null hypothesis that the total number of fires is Normally distributed is not rejected. 

\begin{figure}[H]
	\begin{center}
		\includegraphics[width=3in]{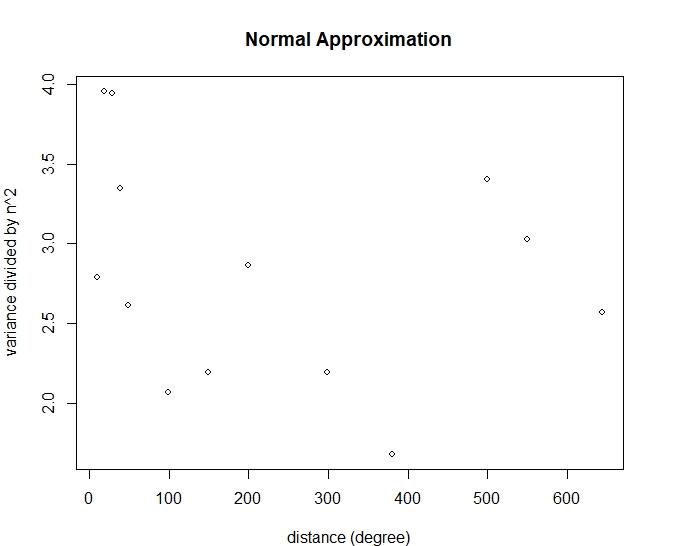}
	\end{center}
	\caption{The variances of fire counts divided by $n^2$}
	\label{fig:varrate}      
\end{figure}

\begin{figure}[H]
	\begin{center}
		\includegraphics[width=3in]{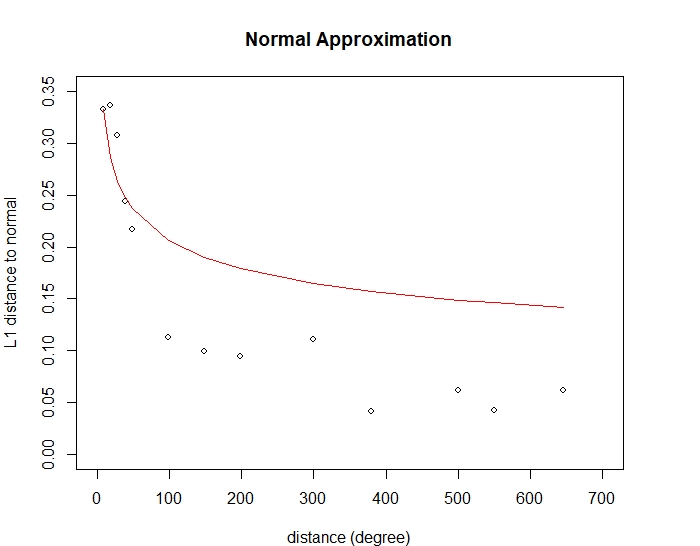}
	\end{center}
	\caption{The $L^1$ distances between the standardized fire counts and a standard Normal. distrbution}
	\label{fig:normapprox}      
\end{figure}

We end this section by checking the variance assumption in assumption \eqref{as2} and compare the simulation’s $L^1$ distance to the bound in \eqref{mainbound2}. Figure \ref{fig:varrate} shows a random pattern for the variances of fire counts divided by $n^2$, which agrees with assumption \eqref{as2}. Figure \ref{fig:normapprox} shows the $L^1$ distances between the standardized fire counts from our simulation and a standard Normal distribution . We also plot the rate of convergence of $n^{d/(4d+2)}$ from Theorem \ref{mainthm2}. The rate from our simulation tends to follow the rate from our main bound. 

\subsection{The effect of changing the covariance decay rate}

The main assumptions that we made in this work is that the fire counts for two distinct areas are positively associated and their covariance decays exponentially with respect to their distance from each other. Nevertheless, the decay rate may vary from country to country or area to area and may also change over time. Therefore, in this section, we are interested in examining the Normal approximation results for different covariance decay rates through the parameter $\lambda_c$ from \eqref{covpat}. 

\begin{figure}[H]
	\begin{center}
		\includegraphics[width=4in]{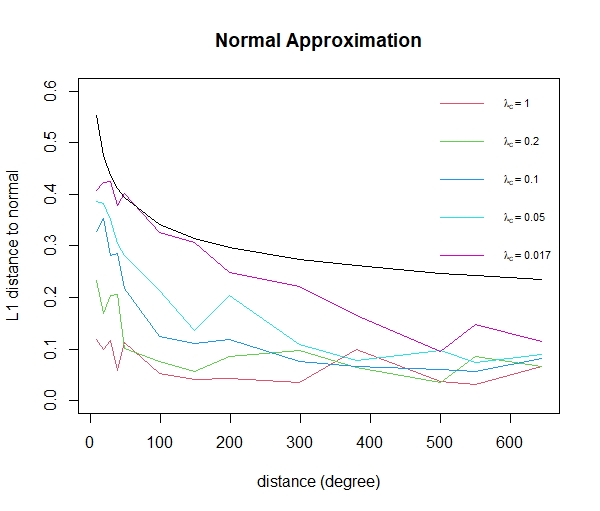}
	\end{center}
	\caption{Normal approximations for different values of $\lambda_c$ }
	\label{fig:normapp}      
\end{figure}

Figure \ref{fig:normapp} shows that the $L^1$ distances between the standardized fire counts from our simulation and a standard Normal distribution for areas ranging in size from 9 to 645 square degrees when $\lambda_c = 1,0.2,0.1,0.05,0.017$. The black line shows the rate of convergence of $n^{d/(4d+2)}$ from Theorem \ref{mainthm2}. Intuitively, the distribution should be closer to the Normal distribution when $\lambda_c$ is larger, as a larger value implies that the fire counts in the areas are closer to independence. This reasoning agrees with the simulation results shown in Figure \ref{fig:normapp}.

\section{Discussion and conclusion} \label{sec:sum}

In this work, we break our contribution into two parts: developing theories and simulating Thailand fire counts. For the first part, we generalize the normal approximation of a functional of the associated point processes in \cite{Wir19} by relaxing the assumption that the covariance decays exponentially  everywhere to exempting local neighborhoods. Then we apply the main result (Theorem \ref{localplus}) to permanental Cox processes, which are known to be positively associated. Unlike the main theorem involving hypercubes of size $n^d$  in $\mathbb{R}^d$, Theorem \ref{mainthm2} is flexible in that it covers Normally approximating a functional of permanental Cox processes in an area of any shape of size $n^d$, which allows us to apply this result to the fire counts in Thailand. 

For the simulation, we use a permanental Cox process to simulate fire counts in Thailand. We assume that the fire counts for two distinct areas are positively correlated. Moreover, we assume that their covariance decays exponentially outside of some small local neighborhood. We use the GISTDA Thailand fire dataset collected by Thailand’s satellites from 2007 to 2020 to estimate the parameters for the process. The $L^1$ distance between the total number of fires from our simulated fire counts and the Normal distribution agrees with the rate from the bound in Theorem \ref{mainthm2}. By varying the covariance parameter in \eqref{covpat}, we can vary the $L^1$ distances from our simulated results. Nevertheless, Normality still holds, and the rates do not contradict the one in our main theorem.

In real-world applications, the new approach proposed in this research may be used to model any natural disaster incidents or model claims by area for property insurance firms. Thus, it benefits policymakers both in government and the private sector in terms of managing risk. Moreover, future researchers can use our idea to simulate natural disaster incidents or insurance claims with some other point processes that are not necessarily well-known but satisfy the assumptions of Theorem \ref{localplus}. 

\section*{Acknowledgment}
The authors would like to thank the Geo-Informatics and Space Technology Development Agency (GISTDA) for use of the fire dataset.

\section*{Declarations}

\begin{flushleft}
	\textbf{Funding:} This work was supported financially by the TSRI Fundamental Fund 2020 (Grant Number: 64A306000047).
	
	\textbf{Conflicts of interest/Competing interests:} No potential conflicts of interest were reported by the authors.
	
	\textbf{Availability of data and material:} This dataset is under the license of the Geo-Informatics and Space Technology Development Agency.
	
	\textbf{Code availability:} Available upon request.
	
	\textbf{IRB approval:} Not applicable.
\end{flushleft}

\def\cprime{$'$}


\begin{thebibliography}{9 }
	
	\bibitem{Godfred14} Abledu, G. K., Dadey, E. and Kobina , A., Probability modeling and simulation of insurance claims in Ghana, 
	{\it Global Journal of Commerce \& Management Perspective}, {\bf 3}(5) (2014) 41--49.
	
	\bibitem{Albrecher2021} Albrecher, H., Araujo-Acuna, J. C. and Beirlant, J., Fitting non-stationary Cox process: an application to fire insurance data, {\it North American Actuarial Journal}, {\bf 0}(0) (2020) 1--28.
	
	\bibitem{Bar88} Barbour, A. D., 
	Stein's method and Poisson process convergence,
	{\it Journal of Applied Probability}, {\bf 25} (1988) 175--184.
	
	\bibitem{BB92} Barbour, A. D. and Brown, T. C., 
	Stein's method and point process approximation,
	{\it Stochastic Processes and their Applications}, {\bf 43}(1) (1992) 9--31.
	
	\bibitem{Bartl2020} B\"{a}rtl, M. and Krummaker, S., 
	Prediction of claims in export credit finance: a comparison of four machine learning techniques,
	{\it Risks}, {\bf 8}(22) (2020) 1--29.
	
	\bibitem{CGS11} Chen, L.Y.H., Goldstein, L. and Shao, Q.M., {\it Normal Approximation by Stein's Method}, Springer, New York, 2011.
	
	\bibitem{ChenShao2004} Chen, L.Y.H. and Shao, Q.M., Normal approximation under local dependence, 
	{\it Annals of Probability}, {\bf 32} (2004) 1985--2028.
	
	\bibitem{CX04} Chen, L.Y.H. and Xia, A., 
	Stein's method, Palm theory and Poisson process approximation,
	{\it Annals of Probability}, {\bf 32}(3) (2004) 2545--2569.
	
	\bibitem{CX06} Chen, L.Y.H. and Xia, A., 
	Poisson process approximation: from Palm theory to Stein's method,
	{\it IMS Lecture Notes-Monograph Series: Time Series and Related Topics}, {\bf 52} (2006) 236--244.
	
	\bibitem{CX11} Chen, L.Y.H. and Xia, A., 
	Poisson process approximation for dependent superposition of point processes,
	{\it Bernoulli}, {\bf 17}(2) (2011) 530--544.
	
	\bibitem{Eis14} Eisenbaum, N., 
	Characterization of positively correlated squared Gaussian processes,
	{\it The Annals of Probability}, {\bf 42}(2) (2014) 559--575.
	
	\bibitem{Adrea18} Gabrielli, A. and W\"{u}thrich, M. V., An individual claims history simulation machine, 
	{\it Risks}, {\bf 6}(29) (2018) 1--33.
	
	\bibitem{Heffernan2007} Heffernan, J.E. and Resnick, S.I., Limit laws for random vectors with an extreme component, {\it Annals of Applied Probability}, {\bf 17} (2007) 537--571.

\bibitem{Jessen2011} Jessen A.H., Mikosch T. and Samorodnitsky G., Prediction of outstanding payments in a Poisson cluster model, {\it Scandinavian Actuarial Journal}, {\bf 3}(2011) 214--237.

\bibitem{Keef2009} Keef, C., Svensson C. and Tawn, J.A., Spatial dependence in extreme river flows and precipitation for Great Britain, {\it Journal of Hydrology}, {\bf 378} (2009) 240--252.

\bibitem{Caroline13} Keef, C., Tawn, J. A. and Lamb, R., Estimating the probability of widespread flood events, {\it Environmetrics}, {\bf 24} (2013) 13--21.

\bibitem{Lamb2005} Lamb, R., Rainfall-runoff modelling for flood frequency estimation, 
{\it Encyclopedia of Hydrological Sciences}, Anderson MG (eds). John Wiley \& Sons(2005) 1923--1954.

\bibitem{LSY20} Last, G., Szekli, R., Yogeshwaran, D., 
Some remarks on associated random fields, random measures and point processes,
{\it ALEA Latin American Journal of Probability and Mathematical Statistics}, {\bf 17} (2020) 355--374.

\bibitem{MM05} McCullagh, P., M\o ller, J., 
The permanental process, {\it Advances in Applied Probability}, {\bf 38}(4) (2006) 873--888.

\bibitem{MY06} McCullagh, P., Yang, J., 
Stochastic classification models,
{\it International Congress of Mathematicians}, {\bf 3} (2006) 669--686.

\bibitem{MESIT18} Mustafa, H. A., Ekti, A. R., Shakir, M. Z., Imran, M. A., Tafazolli, R., 
Intracell interference characterization and cluster interference for D2D communication,
{\it IEEE Transactions on Vehicular Technology}, {\bf 67}(9) (2018) 8536--8548.




\bibitem{PDL17} Poinas, A., Delyon, B. and Lavancier, F., 
Mixing properties and central limit theorem for associated point processes,
{\it Bernoulli,} {\bf 25}(3) (2019) 1724--1754.

\bibitem{Quan2018} Quan, Z. and Valdez, E.A., 
Predictive analytics of insurance claims using multivariate decision trees,
{\it Dependence Modeling}, {\bf 6} (2018) 377--407.

\bibitem{Ross11} Ross, N.
Fundamentals of Stein's method.
{\it Probability Surveys}, {\bf 8} (2011) 210--293.

\bibitem{Schoenberg2007} Schoenberg, F.P., Chang, C., Keeley, J., Pompa, J., Woods, J. and Xu, H., A critical assessment of the burning index in Los Angeles County, California,   
{\it International Journal of Wildland Fire}, {\bf 16}(4) (2007) 473--483.

\bibitem{Stein72} Stein, C., 
A bound for the error in the normal approximation to the distribution of a sum of dependent random variables.
{\it Proceedings of the Sixth Berkeley Symposium on Mathematical Statistics and Probability} University of California Press, {\bf 2} (1972) 210--293.

\bibitem{Gustav2015} T\"{o}rnqvist, G., Modelling insurance claims with spatial point process: an applied case-control study to improve the use of geographical information in insurance pricing, {\it Master thesis, Ume{\r a} University}, 2015


\bibitem{Wir19} Wiroonsri, N., 
Normal approximation for associated point processes via Stein's method with applications to determinantal point processes,
{\it Journal of Mathematical Analysis and Applications}, {\bf 480}(1) (2019) 123396.

\bibitem{Wuthrich2018} W\"{u}thrich M. V., Machine learning in individual claims reserving, 
{\it Scandinavian Actuarial Journal}, (2018) 1--16.

\bibitem{Xu2011} Xu, H. and Schoenberg F.P., Point process modeling of wildfire hazard in Los Angeles County, California,  
{\it The Annals of Applied Statistics}, {\bf 5}(2011) 684--704.

\bibitem{YMM12} Yang, J., Miescke, K., McCullagh, P., 
Classification based on a permanental process with cyclic approximation,
{\it Biometrika}, {\bf 99}(4) (2012) 775--786.

\end{thebibliography}
\end{document}